\input mssymb
\def\proof{\medbreak \noindent{\bf Proof. }}

\def\res {{\rm Res}}
\def\cor {{\rm Cor}}
\def\ind {{\rm Ind}}

\def\qed{\hfill${\scriptstyle\blacksquare}$}

\def\book#1/#2/#3/#4/#5/#6/ {\item{[#1]} #2{\sevenrm #3}. {\it #4}.\
#5, (#6).\par\smallskip}
\def\paper#1/#2/#3/#4/#5/#6/#7/#8/ {\item{[#1]} #2{\sevenrm #3}. #4.
{\it #5} {\bf #6} (#7), #8.\par\smallskip}
\def\prepaper#1/#2/#3/#4/#5/ {\item{[#1]} #2{\sevenrm #3}. #4.
#5\par\smallskip}
\baselineskip= 20 pt

\def\ch{{\rm Ch}}
\def\gee{G}

\def\ggen#1#2{{\gee(#1,#2)}}
\def\gne{{\ggen{n}{\epsilon}}}
\def\bgen#1{{{\rm B}#1}}
\def\bgee{\bgen \gee}

\def\co#1#2#3{{H^{#1}(#2;#3)}}
\def\coz#1#2{\co{#1}{#2}{\Bbb Z}}

\def\coev#1{\coz{{\rm even}}{#1}}

\def\ala{1}
\def\alb{2}
\def\atiyah{3}
\def\blackburn{4}
\def\brown{5}
\def\ccsand{6} 
\def\land{7}
\def\leary{8}
\def\lewis{9}
\def\tyone{10}
\def\tytwo{11}
\def\meta{12}
\def\riemann{13}
\def\admiss{14}
\def\thoc{15}
\def\weiss{16}

\magnification = \magstep 1

\centerline{{\bf Some examples in the integral and Brown-Peterson cohomology of
$p$-groups.}}
$$\hbox{
\vbox{\hbox{I. J. Leary}
\hbox{School of Mathematical Sciences,}
\hbox{Queen Mary and Westfield College,}
\hbox{London.}
\hbox{E1 4NS}}
\qquad\vbox{\hbox{N. Yagita}
\hbox{Musashi Institute of Technology,}
\hbox{Tamazatumi,}
\hbox{Setagaya,}
\hbox{Tokyo.}}}$$

\beginsection{Introduction.}\par
For a finite group $G$, we define the Chern ring, $\ch(G)$, to be the subring
of $\coev G$ generated by Chern classes of representations of $G$. We say that
$G$ has $p$-rank $n$ if $n$ is maximal such that $G$ contains a subgroup
isomorphic to $(C_p)^n$. In [\atiyah] Atiyah showed that for any finite group
$G$, $K^0(\bgee)$ is the completion of the representation ring of $G$ with
respect to a certain topology. The filtration of $K^0(\bgee)$ given by the
$E_\infty$ page of the Atiyah-Hirzebruch spectral sequence (AHSS) gives rise to
a filtration of the representation ring of $G$. Atiyah conjectured that this
filtration coincided with another filtration defined algebraically, and
remarked that this conjecture is equivalent to the conjecture that $\ch(G)$
maps onto the $E_\infty$ page of the AHSS. (It is clear that $\ch(G)$ consists
of universal cycles because the AHSS for $\bgen{U(n)}$ collapses.) Weiss
discovered that the alternating group $A_4$ gives a counterexample to this
conjecture [\weiss], and Thomas has exhibited many counterexamples,
all of which have order divisible by more than one prime [\admiss].
\par
Thomas showed that the split metacyclic $p$-groups and various other $p$-groups
of $p$-rank two have the property that the Chern subring is the whole of the
even degree integral cohomology, and conjectured that this property would hold
for all $p$-groups of $p$-rank two [\meta], [\riemann], [\thoc].
The group $A_4$ shows that the conjecture cannot be extended to groups of
non-prime power order. AlZubaidy claimed
to have verified this conjecture, but some of his proofs are flawed [\ala],
[\alb]. Recently Huebschmann and Tezuka-Yagita have shown that $\ch(G)=
\coev\gee$ for any metacyclic $p$-group $G$ [\ccsand], [\tytwo].
For $p\geq 5$ Blackburn's classification [\blackburn] implies that the
only $p$-groups of $p$-rank two not covered by the above theorems are the
groups $\gne$ defined below. We show that $\ch(\gne)$ is strictly contained in
$\coev\gne$ for each such group.
These groups also afford $p$-group counterexamples to the conjecture of Atiyah
described above.
\par
Similar calculations may be made in the Brown-Peterson cohomology rings of
these groups. These enable us to give a negative answer to a question of
Landweber [\land], who asked if Chern classes generate the Brown-Peterson
cohomology of every $p$-group.
\par
\beginsection{The examples.}
\par
The groups which we shall consider may be presented as
$$\gne=\langle A,B,C\mid A^p=B^p=C^{p^{n-2}}=[B,C]=1\quad [A,C^{-1}]=B\quad
[B,A]=C^{\epsilon p^{n-3}}\rangle,$$
where $p$ is a prime not equal to 2 or 3, $n\geq 4$, and
for fixed $p$ and $n$ there are two isomorphism classes of such groups,
depending whether $\epsilon$ is either 1
or a quadratic non-residue modulo $p$. The group $\gne$ has order $p^n$.
In the sequel we shall refer to $\gee$ instead of $\gne$ unless the values
of $n$ and $\epsilon$ are important. The subgroup $M$ generated by $B$ and
$C$ is maximal (and hence normal) and is isomorphic to $C_p\oplus C_{p^{n-2}}$.
We define one dimensional representations $\theta$ and $\phi$ of the group $M$
by
$$\eqalign{\theta : B^jC^k&\mapsto \exp({2\pi ij/p})\cr
\phi : B^jC^k&\mapsto \exp({2\pi ik/{p^{n-2}}}).}$$
The action of the quotient group $\gee/M$ on the representation ring of $M$ is
that conjugation by $A$ sends $\theta$ to $\theta\otimes\phi^{\otimes
p^{n-3}}$ and sends $\phi$ to $\phi\otimes\theta^{\otimes\epsilon}$. Later we
shall define elements of $\coz 2 M$ and $BP^2({\rm B}M)$ as Chern classes, and
the action of $\gee/M$ on these elements will be determined by its action on the
representations $\theta$ and $\phi$.
\par
The group $\gee$ has only 1- and $p$-dimensional irreducible representations
because it has an abelian subgroup ($M$ in fact) of index $p$. A one
dimensional representation of $\gee$ must restrict trivially to $\langle B
\rangle$, and a $p$-dimensional representation of $\gee$ restricts to
$\langle B \rangle$ as either $p$ copies of the same representation of
$\langle B \rangle$, or as the sum of one copy of each of the one-dimensional
representations of $\langle B \rangle$. The examples $\ind^\gee_M(\theta)$
and $\ind^\gee_M(\phi)$ show that both these alternatives do occur.
\par
Now define generators $\beta$, $\gamma$ for $\coz 2 M$ by
$$\beta= c_1(\theta)\qquad \gamma= c_1(\phi),$$
$$\hbox{so that} \qquad \coev M \cong \Bbb Z[\beta,\gamma]/(p\beta,
p^{n-2}\gamma),$$
and let $\beta'$ be the restriction to $\langle B \rangle$ of $\beta$, so that
$$\coz * {\langle B\rangle}\cong\Bbb Z[\beta']/(p\beta').$$
\par
\proclaim Lemma 1. With notation as above, the image of $\ch(\gee)$ under
restriction to $\langle B\rangle$ is the subring of $\coz *
{\langle B \rangle}$ generated by $\beta'^{p-1}$ and $\beta'^p$. For all
$m \geq 0$,
$$ \beta'^{m+p-1}=-\res^\gee_{\langle B \rangle}\cor^\gee_M(\gamma^{p-1}
\beta^m).$$
\par
\proof If $\rho$ is a 1-dimensional representation of $\gee$ then its Chern
class restricts trivially to $\langle B\rangle$. If $\rho$ is a
$p$-dimensional representation of $\gee$, then either $\rho$
\def\be{{\langle B\rangle}}
restricts to $\be$ as $p$-copies of the same representation, in which case
$$\res(c.(\rho))=(1+i\beta')^p=1+i\beta'^p,$$
or as one copy of each representation, in which case
$$\res(c.(\rho))=\prod^{p-1}_{i=0}(1+i\beta')=1-\beta'^{p-1}.$$
By applying the double coset formula we see that
$$\eqalignno{\res^\gee_\be\cor^\gee_M(\gamma^{p-1}\beta^m)&=
\res^M_\be(\sum^{p-1}_{i=0}c_{A^i}^*(\gamma^{p-1}\beta^m)\cr
&=\res^M_\be(\sum_{i=0}^{p-1}(\gamma+i\epsilon\beta)^{p-1}(\beta+ip^{n-3}
\gamma)^m)\cr
&=-\beta'^{m+p-1}&{\scriptstyle \blacksquare}}$$
\par
{\bf Remarks.} The image of $\res^\gee_\be$ is precisely the subring of
$\coz * \be $ generated by $\beta'^{p-1}$, $\beta'^p$,
$\beta'^{p+1},\ldots,\beta'^{2p-3}$. One way to show this is by considering
the subgroup $N$ of $\gee$ generated by $A$ and $B$. This subgroup is normal
in $\gee$, and is the non-abelian group of order $p^3$ and exponent $p$.
Using Lewis' calculation of $\coz * N$ [\lewis], it may be shown that the image
of $\coz * N^{\gee/N}$ under restriction to $\be$ does not contain $\beta'^i$
for $i < p-1$.
\proclaim Corollary 2. $\ch(\gee)$ is strictly contained in $\coev \gee$.
Moreover, $\ch(\gee)$ does not map onto $\coz * \gee$ modulo its nilradical.
\par
\proof We know that $\beta'^{p+1}$ is in $\res^\gee_\be(\coev \gee)$, but not
in $\res^\gee_\be(\ch(\gee))$. \qed
\par
\proclaim Corollary 3. In the AHSS for $\gee$, write $B_\infty(\gee)$ for the
universal boundaries, and $Z_\infty(\gee)$ for the universal cycles. Then
$\ch(\gee)+B_\infty(\gee)$ is strictly contained in $Z_\infty(\gee)$.
\par
\proof The AHSS for $\be$ collapses, so $\res^\gee_\be(B_\infty(\gee))$ is
trivial. Corestrictions of Chern classes are universal cycles, so
$\beta'^{p+1}\in \res^\gee_\be(Z_\infty(\gee))$, but
$\beta'^{p+1}\notin \res^\gee_\be(\ch(\gee)+B_\infty(\gee))$.\qed\par
For any generalised cohomology theory $\cal H$ and any group $K$, we may
define $\ch_{\cal H}(K)$ to be the subring of
${\cal H}^*({\rm B}K)$ generated by
$\rho^*({\cal H}^*({\rm B}U))$ for all representations $\rho$ of $K$ in a
unitary group $U$.
We now give a result analogous to Corollary~2 for Brown-Peterson cohomology.
\proclaim Lemma 4. $\ch_{BP}(\gee)$ is strictly contained in $BP^*(\bgee)$.
\par
\proof As in the integral cohomology case, define elements $\beta$ and
$\gamma$ in $BP^2(M)$ by $\beta=c_1(\theta)$, $\gamma=c_1(\phi)$, and also
define $\beta'=\res^M_\be(\beta)$, so that
$$BP^*({\rm B}M)\cong BP_*[[\beta,\gamma]]/([p]\beta,[p^{n-2}]\gamma),
\qquad BP^*({\rm B}\be)\cong BP_*[[\beta']]/([p]\beta'),$$
where $[r]x$ stands for the $BP$ formal group sum of $r$ copies of $x$.
Let `$\equiv$' stand for congruence modulo the ideal of $BP^*({\rm B}\be)$
generated by $p$, $v_1$, $v_2,\ldots$. As in Lemma~1, if $\rho$ is a
$p$-dimensional representation of $\gee$, then either
$$\eqalign{\res^\gee_\be(c.(\rho))&=(1+[i]\beta')^p\equiv
1+i\beta'^p\cr
\hbox{or}\quad\res^\gee_\be(c.(\rho))&=\prod_{i=0}^{p-1}
(1+[i]\beta')\equiv 1-\beta'^{p-1}.}$$
Also, we have that
$$\eqalignno{\res^\gee_\be\cor^\gee_M(\gamma^{p-1}\beta^2)&=
\res^M_\be(\sum_{i=0}^{p-1}(\gamma+_{BP}[i\epsilon]\beta)^{p-1}
(\beta+_{BP}[ip^{n-3}]\gamma)^2\cr
&\equiv-\beta'^{p+1}.&{\scriptstyle \blacksquare}}$$
\par
Our original proofs of these results involved calculation with $BP^*({\rm B}N)$,
which has been determined by Tezuka-Yagita [\tyone], and with the integral
cohomology of the non-abelian maximal subgroups of $\gee$, determined by
Leary [\leary]. Using these methods we obtain more information concerning
$BP^*(\bgee)$ and $\coz * \gee$, which we intend to publish later.
\par
We thank Dr.~C.~B.~Thomas for suggesting investigations of the Chern
subring and for his encouragement. The second author also thanks M.~Tezuka
and M.~Kaneda for many suggestions. The work of the first author has been
supported by the SERC.
\par \bigskip \noindent
{{\rm REFERENCES}}\par \medskip
\paper \ala/K.~A/L{\rm Z}UBAIDY/Metacyclic $p$-groups and Chern classes/
Illinois J. Math./26/1982/423--431/
\paper \alb/K.~A/L{\rm Z}UBAIDY/Rank 2 $p$-groups, $p>3$, and Chern classes/
Pacific J. Math./103/1982/259--267/
\paper \atiyah/M.~F.~A/TIYAH/Characters and the cohomology of finite groups/
Publ. Math. IHES/9/1961/23--64/
\paper \blackburn/N.~B/LACKBURN/Generalisations of certain elementary
theorems on $p$-groups/Proc. London Math. Soc./11/1961/1--22/
\book \brown/K.~S.~B/ROWN/Cohomology of groups/Springer Verlag/1982/
\prepaper \ccsand/J.~H/UEBSCHMANN/Chern classes for metacyclic groups/preprint
(1988)./
\paper \land/P.~L/ANDWEBER/Elements of infinite filtration in complex
cobordism/Math. Scand./30/1972/223--6/
\prepaper \leary/I.~J.~L/EARY/The integral cohomology rings of some
$p$-groups/Submitted./
\paper \lewis/G.~L/EWIS/The integral cohomology rings of groups of order $p^3$/
Trans. Amer. Math. Soc./132/1968/501--529/
\paper \tyone/M.~T/EZUKA {\tenrm and N.~Y}AGITA/Cohomology of finite groups and
Brown-Peterson cohomology/Springer LNM/1370/1989/396--408/
\paper \tytwo/M.~T/EZUKA {\tenrm and N.~Y}AGITA/Cohomology of finite groups
and Brown-Peterson cohomology II/Springer LNM/1418/1990/57--69/
\paper \meta/C.~B.~T/HOMAS/Chern classes and metacyclic $p$-groups/
Mathematika/18/1971/196--200/
\paper \riemann/C.~B.~T/HOMAS/Riemann-Roch formulae for group representations/
Mathematika/20/1973/253--262/
\paper \admiss/C.~B.~T/HOMAS/Filtrations on the representation ring of a finite
group/Cont. Math./19/1983/397--405/
\book \thoc/C.~B.~T/HOMAS/Characteristic classes and the cohomology of finite
groups/Cambridge University Press/1986/
\prepaper \weiss/E.~A.~W/EISS/Bonner Mathematische Schriften/1969/

\end